\documentclass[11pt]{amsart}
\usepackage[english]{babel}
\usepackage[utf8]{inputenc}
\usepackage{amsmath,amssymb,mathtools}
\usepackage{enumitem}
\usepackage{xcolor,tikz-cd}
\usepackage[colorlinks=true,linkcolor=blue,citecolor=blue,urlcolor=magenta,pagebackref=true]{hyperref}

\numberwithin{equation}{section}

\newtheorem{thm}{Theorem}[section]
\newtheorem{prop}[thm]{Proposition}
\newtheorem{lemma}[thm]{Lemma}
\newtheorem{cor}[thm]{Corollary}
\newtheorem{question}[thm]{Question}
\theoremstyle{definition}
\newtheorem{rmk}[thm]{Remark}

\newcommand{\Z}{\mathbb{Z}}
\newcommand{\Q}{\mathbb{Q}}
\newcommand{\F}{\mathbb{F}}

\newcommand{\m}{\mathfrak m}
\newcommand{\GL}{\operatorname{GL}}
\newcommand{\End}{\operatorname{End}}
\newcommand{\Hom}{\operatorname{Hom}}
\newcommand{\tr}{\operatorname{tr}}
\newcommand{\Tr}{\operatorname{Tr}}
\newcommand{\Ext}{\operatorname{Ext}}
\newcommand{\ad}{\operatorname{ad}}
\newcommand{\Sym}{\operatorname{Sym}}

\renewcommand{\phi}{\varphi}

\DeclareFontFamily{U}{wncy}{}
\DeclareFontShape{U}{wncy}{m}{n}{<->wncyr10}{}
\DeclareSymbolFont{mcy}{U}{wncy}{m}{n}
\DeclareMathSymbol{\Sh}{\mathord}{mcy}{"58} 

\title[The lifting problem for Galois representations of global fields]{The lifting problem for Galois representations of global fields}
\author{Federico Scavia}
\address{CNRS, Institut Galil\'{e}e, Universit\'{e} Sorbonne Paris Nord,
	99 avenue Jean-Baptiste Cl\'{e}ment, 93430 Villetaneuse, France}
\email{scavia@math.univ-paris13.fr}
\date{October 2026}

\begin{document}
	
	\begin{abstract}
		We prove that all mod $2$ Galois representations of a global field lift modulo $4$. In contrast, for every odd prime $p$, we exhibit examples of global Galois representations modulo $p$ that do not lift modulo $p^2$. Our methods also allow us to prove liftability of generalized quaternion groups and hence to complete the classification of liftable finite groups initiated by Khare and Merkurjev.
	\end{abstract}
	
	\maketitle

	\section{Introduction}
	
	Let $G$ be an abstract or profinite group. If $G$ is profinite, we assume that every homomorphism $G\to H$, where $H$ is an abstract group, is continuous once we endow $H$ with the discrete topology, and that every $G$-module is discrete. Following Khare--Larsen \cite{KL24} and Khare--Merkurjev \cite{KM26}, for a prime $p$ we say that $G$ is
	
	\begin{itemize}
		\item $p$-\emph{liftable} if for every finite field $\F_q$ of characteristic $p$ and every integer $n\geq 1$, every homomorphism $G\to \GL_n(\F_q)$ lifts to $\GL_n(W_2(\F_q))$;
		\item \emph{liftable in
			characteristic $p$} if for every commutative local ring $R$ with square-zero maximal ideal $\m$ and residue field $k=R/\m$ of characteristic $p$ and every integer $n\geq 1$, every homomorphism $G\to\GL_n(k)$ lifts to $\GL_n(R)$;
		\item  \emph{liftable} if it is liftable in characteristic $p$, for every prime $p$.
	\end{itemize} 
	
	Khare and Serre \cite{Kha97} proved that, for every finite field $\F_q$, every homomorphism $G\to \GL_2(\F_q)$ lifts to $\GL_2(W_2(\F_q))$: more precisely, Khare proved this when $F$ is a number field, and Serre observed that Khare's argument only involves Kummer theory and so applies to an arbitrary field \cite[Remark~2, p.~392]{Kha97}. They then asked whether absolute Galois groups are $p$-liftable for every prime $p$.  De Clercq and Florence \cite{DCF22} later extended the theorem of Khare and Serre to arbitrary coefficient fields of characteristic $p$, and also showed that, for $n\in \{3,4\}$, all $n$-dimensional $\F_2$-representations of an absolute Galois group lift to $\Z/4\Z$. For short alternative proofs of these results, see \cite[\S3]{MSMA26}.
	
	The answer to the question of Khare--Serre turns out to be negative in general. Indeed, for every odd prime $p$, every $n\geq3$, and every field $F$ containing $\mu_p$, Merkurjev and the author constructed $n$-dimensional an representation of the absolute Galois group of $F(x_1,\ldots,x_p)$ over $\F_p$ which does not lift to $\Z/p^2\Z$, and they also constructed counterexamples in dimension $\geq 5$ over $\F_2$; see \cite[Theorem~1.4]{MSJAMS26} and \cite[Theorem~1.1]{MSMA26}. The dimension of the representations considered in these counterexamples are optimal for all primes $p$, in view of the aforementioned results of Khare--Serre and De Clercq--Florence. However, the following question had so far remained open.
	
	\begin{question}\label{question:global}
		Are absolute Galois group of global fields liftable?
	\end{question}
	
	Question~\ref{question:global} belongs to the same circle of ideas as the deformation theory of Galois representations and its applications to modularity, notably in the work of Ramakrishna \cite{Ram99,Ram02} and Khare--Wintenberger \cite{KW09a,KW09b}. However, an important difference is that here we do not impose additional condition on the lifts. Another closely related lifting problem arises in the Massey Vanishing Conjecture of Min\'a\v{c}--T\^an \cite{MT17}, which concerns lifting Galois representations through central extensions of $p$-groups related to the group of upper unitriangular matrices over $\F_p$. This conjecture was proved for number fields by
	Harpaz--Wittenberg \cite{HW23}, but remains open over arbitrary fields, although several cases have been proved; see \cite{MT16,EM17,MSJEMS26,MSENS25,MSPLMS25,MSComp25}.
	
	By a theorem of B\"ockle, every continuous mod $p$ representation of the absolute
	Galois group of a local field lifts modulo $p^2$; see
	\cite[Theorem~1.1]{Boc03}, and see \cite{CDF26} for an extension to
	complete flags. B\"ockle also gave sufficient conditions for lifting
	Galois representations global fields \cite{Boc03}, while Khare and Larsen proved that
	every three-dimensional mod $p$ representation of the absolute Galois
	group of a number field containing $\mu_{p^2}$ lifts modulo $p^2$;
	see \cite[Theorem~5.4]{KL24}.

	Our first main theorem is the following.
	
	\begin{thm}\label{thm:main-even}
		The absolute Galois group of every global field is liftable in characteristic $2$.
	\end{thm}
	
	As the next result shows, Theorem~\ref{thm:main-even} cannot be extended to any odd prime $p$, and hence in particular Question~\ref{question:global} has a negative answer in general.
	
	\begin{thm}\label{thm:main-odd}
		Let $p$ be an odd prime, let $F$ be a global field of characteristic
		different from $p$ such that $\mu_{p^2}\not\subset F(\mu_p)$, and let
		$m\coloneqq[F(\mu_p):F]$. Then there exists a continuous homomorphism
		$\rho\colon G_F\to\GL_{2mp}(\F_p)$ such that for every place $v$ of $F$
		the restriction $\rho_v$ lifts to $\GL_{2mp}(\Z_p)$, whereas $\rho$ does
		not lift to $\GL_{2mp}(\Z/p^2\Z)$.
	\end{thm}
	
	For a number field $F$, the hypothesis of Theorem~\ref{thm:main-odd}
	holds whenever $p$ does not divide $[F:\Q]$. Indeed, if
	$\mu_{p^2}\subset F(\mu_p)$, then
	\[
	p(p-1)\mid[F(\mu_p):\Q]
	=[F(\mu_p):F][F:\Q]\mid(p-1)[F:\Q],
	\]
	so $p\mid[F:\Q]$. For $F=\Q$, Theorem~\ref{thm:main-odd} gives, for every odd $p$, an $\F_p$-linear representation of dimension $2p(p-1)$ which does not lift to $\Z/p^2\Z$: thus $G_{\mathbb Q}$ is liftable in characteristic $2$ but it is not liftable in characteristic $p$ for any odd prime $p$.
	
	\subsection{Liftable finite groups}
	
	Our methods also apply to the classification of finite liftable groups, a problem which was first posed by Fakhruddin \cite[Question~1.2]{KL24}.
	
	\begin{question}[Fakhruddin]\label{qu:fakhruddin}
		Which finite groups are liftable?
	\end{question}
	
	Khare and Merkurjev proposed a
	complete classification \cite[Theorem~1.1]{KM26v1}. Their argument
	included the assertion that some representation of the quaternion group
	$Q_8$ over $\F_2$ does not lift to $\Z/4\Z$;
	see \cite[Proposition~4.2]{KM26v1}. We prove that the opposite holds,
	and more generally that every generalized quaternion group is liftable. More precisely, let $m\geq2$ be a power of $2$, and set
	\[
	Q_{4m}\coloneqq\langle a,b\mid a^{2m}=1,\ b^2=a^m,\ bab^{-1}=a^{-1}\rangle,
	\]
	the generalized quaternion group of order $4m$.
	
	\begin{thm}\label{thm:quaternion}
		Let $m\geq2$ be a power of $2$. Then $Q_{4m}$ is liftable.
	\end{thm}
	
	Combining this result with the work of Khare and Merkurjev gives a
	complete answer to Question~\ref{qu:fakhruddin}, amending the
	classification proposed in \cite[Theorem~1.1]{KM26v1}.
	
	\begin{cor}\label{cor:classification}
		Let $G$ be a finite group. Then $G$ is liftable if and only if the
		following two conditions hold:
		\begin{enumerate}
			\item the order of $G$ is equal to $2^r$ or to $3\cdot2^r$ for some
			integer $r\geq0$;
			\item the Sylow $2$-subgroups of $G$ are cyclic or generalized quaternion.
		\end{enumerate}
	\end{cor}
	
	\subsection{Proof outlines}
	The proofs of Theorems~\ref{thm:main-even} and \ref{thm:quaternion} rely on a common principle. For a field $k$ of characteristic $2$ and a finite-dimensional $G$-representation $V$ over $k$, the obstruction to lifting $V$ to $W_2(k)$ is a class
	\[
	o(V)\in H^2(G,\End(V)^{(1)}),
	\]
	where $\End(V)^{(1)}$ is the Frobenius twist of the endomorphism module;
	see \eqref{eq:obstruction-pullback}. We observe that $o(V)$ can be
	described using the second coefficient of the characteristic polynomial on $\mathrm{End}(V)$. 
	More precisely, for every finite-dimensional
	$G$-representation $W$, there is a natural exact sequence
	\begin{equation}\label{eq:quadratic-sequence-intro}
		0\longrightarrow W^{(1)}
		\longrightarrow \Sym^2(W)
		\longrightarrow W\otimes W
		\longrightarrow \Sym^2(W)
		\longrightarrow 0;
	\end{equation}
	see \eqref{eq:quadratic-sequences}.
	Therefore, every $q\in(\Sym^2 W)^G$ gives rise to a class
	\[
	\kappa_q\in H^2(G,W^{(1)}).
	\]
	We show, by reduction to the universal case, that
	$o(V)=\kappa_{c_2}$, where $c_2$ is the second coefficient
	of the characteristic polynomial on $\End(V)$,
	and where we have identified $\End(V)^\vee$ with $\End(V)$
	via the trace pairing; see Proposition~\ref{prop:comparison}.

	Theorems~\ref{thm:quaternion} and \ref{thm:main-even} will follow from two different properties of the classes $\kappa_q$. First, the pairing
	\[
	H^2(G,W^{(1)})\times H^1(G,(W^{(1)})^\vee)
	\longrightarrow H^3(G,k),
	\qquad
	(\eta,\xi)\longmapsto\langle\eta\cup\xi\rangle,
	\]
	induced by cup product and evaluation, vanishes for
	$\eta=\kappa_q$; see Proposition~\ref{prop:annihilation}. This implies the following result.

	\begin{thm}\label{thm:universal}
		Let $G$ be an abstract or profinite group, let $k$ be a field of
		characteristic $2$, and let $V$ be a finite-dimensional
		$G$-representation over $k$. Then, for every
		$\xi\in H^1(G,\End(V)^{(1)})$, we have
		\[
		\tr_*(o(V)\cup\xi)=0\quad\text{in }H^3(G,k).
		\]
	\end{thm}
		
	For generalized quaternion groups, the cup-product pairing in Theorem~\ref{thm:universal} is perfect by Lemma~\ref{lem:duality}. It follows that $o(V)=0$ for every representation $V$, proving Theorem~\ref{thm:quaternion}. The same argument also gives lifting results for absolute Galois groups of certain two-dimensional local fields and for fundamental groups
	of closed $3$-manifolds; see Section~\ref{sec:ordinary-duality}.
	
	The second property of $\kappa_q$ that we use is that, given an extension of $G$-representations
	\[0\longrightarrow k\longrightarrow E\longrightarrow W\longrightarrow0,\]
	the class $\kappa_q$ admits a canonical lift to $H^2(G,E^{(1)})$ which becomes trivial after restriction to any subgroup where both $\kappa_q$ vanishes and
	the extension splits; see Lemma~\ref{lem:lifting-kappa} for the precise statement.
	Over global fields, this compatibility and Poitou--Tate duality allow us to show in Theorem~\ref{thm:two-S} that the profinite groups $G_{F,S}$, where $S$ is a finite set of places containing all archimedean places and all places above $2$, are liftable in characteristic $2$, and hence to prove Theorem~\ref{thm:main-even}.
	
	We conclude this introduction by outlining the proof strategy of Theorem~\ref{thm:main-odd}. By Lemma~\ref{lem:induction}, we may assume that $\mu_p\subset F$. Starting from two standard extensions of modules over $\F_p[C_p]$, we construct a certain $2p$-dimensional representation $V$ of the bicyclic group $C_p\times C_p$, equipped with an endomorphism $d$ such that $d^2=0$. It is not difficult to show that $V$ lifts to $\Z_p$ everywhere locally. 
	
	For any characters $a,b\colon G_F\to C_p$, the obstruction to lifting the resulting $G_F$-representation modulo $p^2$ is
	\[
	o(V)=d_*\beta(b)\in H^2(F,\End_{\F_p}(V))\simeq\Ext^2_{G_F}(V,V).
	\]
	where $\beta\colon H^1(F,\F_p)\to H^2(F,\F_p)$ denotes the Bockstein homomorphism, and $d_*$ is induced by the $G_F$-equivariant map $\F_p\to\End_{\F_p}(V)$ sending $1$ to $d$. Our goal is to prove that $d_*\beta(b)\neq 0$.
	
	A first obstruction to the vanishing of $o(V)$ can be obtained by considering the trace map $H^2(F,\End_{\F_p}(V))\to H^2(F,\F_p)$; however, since $\tr(d)=0$, this obstruction is trivial on $o(V)$. In order to refine this obstruction, we construct an extension
	\begin{equation}\label{eq:v-p-v}
		0 \longrightarrow V \xlongrightarrow{i} P \xlongrightarrow{\pi} V \longrightarrow 0
	\end{equation}
	and a $(C_p\times C_p)$-equivariant map $\Phi\colon V\to P$ satisfying $\pi\Phi=d$. 
	The point is that, if $o(V)=0$, the class $\Phi_*\beta(b)\in\Ext^2_{G_F}(V,P)$ lies in the kernel of $\pi_*$ and hence, by the exactness of
	\[
	\Ext^2_{G_F}(V,V)\xlongrightarrow{i_*}\Ext^2_{G_F}(V,P)\xlongrightarrow{\pi_*}\Ext^2_{G_F}(V,V),
	\]
	there exists $\Gamma\in \Ext^2_{G_F}(V,V)$ such that $i_*\Gamma=\Phi_*\beta(b)$, and taking the trace of $\Gamma$ gives a class $\gamma\in H^2(F,\F_p)$. 
	
	We construct \eqref{eq:v-p-v} and $\Phi$ so that the extension splits over every order-$p$ subgroup of $C_p\times C_p$ and the traces of its equivariant retractions take prescribed values; see Lemma~\ref{lem:local-traces}. We then choose $a,b$ so that $a_v$ and $b_v$ are linearly dependent at every place $v$. Thus \eqref{eq:v-p-v} splits over $G_{F_v}$, and for every equivariant retraction $r_v\colon P\to V$ of $i$, we have $\gamma_v=\tr(r_v\Phi)\beta(b_v)$, so that Brauer reciprocity for $\gamma$ gives
	\[
	\sum_v\tr(r_v\Phi)\operatorname{inv}_v(\beta(b_v))=0.
	\]
	For our choice of $a$, $b$ and $P$, the traces are zero when $a_v\neq 0$ and equal to $2$ when $a_v=0$ and $b_v\neq 0$. Using Brauer reciprocity for $\beta(b)$, we can find $a$ and $b$ such that the preceding sum equals $2\operatorname{inv}_{v_0}(\beta(b_{v_0}))\neq 0$ for some place $v_0$; see Lemma~\ref{lem:arithmetic}. This contradiction proves that $o(V)\neq 0$.

	\subsection*{Notation}
	
	Let $G$ be an abstract or profinite group, let $k$ be a field, and let $V$ be a finite-dimensional $k$-representation of $G$, corresponding to a homomorphism $\rho\colon G\to\GL(V)$. We write $\ad(\rho)$ for the adjoint representation of $\rho$, that is, $\ad(\rho)\coloneqq\End_k(V)$ where $G$ acts on $\End_k(V)$ by conjugation. If $k$ has characteristic $p>0$, we let $V^{(1)}\coloneqq V\otimes_{k,\mathrm{Fr}}k$ denote the Frobenius twist of $V$, where $\mathrm{Fr}(a)=a^p$. The Frobenius map $V\to V^{(1)}$ sends $v$ to $v\otimes1$, and we endow $V^{(1)}$ with the Frobenius twist of the $G$-action. 
	
	Let $A$ be a $G$-module, which we assume to be discrete if $G$ is profinite.  For every $i\geq0$, we write $Z^i(G,A)$ and $H^i(G,A)$ for the groups of $i$-cocycles and the cohomology groups of degree $i$, respectively, where we use continuous cochains if $G$ is profinite. As mentioned in the Introduction, if $G$ is profinite, we always assume that every homomorphism $G\to H$, where $H$ is an abstract group, is continuous once we endow $H$ with the discrete topology, and that every $G$-module is discrete.
	
	Let $F$ be a field, and let $G_F$ be the absolute Galois group of $F$. For every positive integer $n$ invertible in $F$, we let $\mu_n$ denote the $G_F$-module of $n$th roots of unity. If $F$ is a global field, for every place $v$ of $F$ we let $F_v$ denote the completion of $F$ at $v$. For every object $\alpha$ defined over $F$ (such as a Galois representation, a cochain, or a cohomology class) we write $\alpha_v$ for its restriction to $F_v$.
	
	Let $\mathcal G$ be an affine group scheme over a field $k$, and let $M$ be a representation of $\mathcal G$. We write $H^*_{\mathrm{rat}}(\mathcal G,M)$ for the rational (Hochschild) cohomology of $\mathcal G$ with coefficients in $M$, which is computed by the Hochschild complex, whose $i$-cochains are the morphisms of $k$-schemes $\mathcal G^i\to M$, where $M$ is viewed as a vector group over $k$; see \cite[I.4.14--4.16]{Jan03}. In particular, if $G$ is an abstract or profinite group and $\rho\colon G\to\mathcal G(k)$ is a continuous homomorphism, then precomposition with $\rho$ induces a map $\rho^*\colon H^i_{\mathrm{rat}}(\mathcal G,M)\to H^i(G,M)$, where $G$ acts on $M$ through $\rho$. We equip $M^{(1)}$ with
	the Frobenius twist of the action of $\mathcal G$ on $M$: if the $\mathcal G$-action on $M$ is given by a representation $\rho\colon\mathcal G\to\GL(M)$, the $\mathcal G$-action on $M^{(1)}$ is given by the composite $\mathcal G\xrightarrow{\rho}\GL(M)\xrightarrow{\mathrm{Fr}_*}\GL(M^{(1)})$, where $\mathrm{Fr}_*$ is induced by $\mathrm{Fr}\colon M\to M^{(1)}$.

	\section{Preliminaries on the lifting problem}\label{sec:prelim}
	
	\subsection{The lifting obstruction}
	
	Let $R$ be a commutative local ring with maximal ideal $\m$ such that $\m^2=0$, let $k\coloneqq R/\m$, and let $n\geq1$. Since $\m^2=0$, the ideal $\m$ is a $k$-vector space, and reduction gives a group extension
	\begin{equation}\label{eq:square-zero}
		0\longrightarrow M_n(\m)\longrightarrow\GL_n(R)\longrightarrow\GL_n(k)\longrightarrow1,
	\end{equation}
	where the first map sends $x$ to $I+x$, and the induced $\GL_n(k)$-action on $M_n(\m)$ is given by conjugation.
	For every abstract or profinite group $G$ and every (continuous) homomorphism $\rho\colon G\to\GL_n(k)$, let $o_R(\rho)\in H^2(G,M_n(\m))$ be the class of the pullback of \eqref{eq:square-zero} along $\rho$. By definition, $o_R(\rho)=0$ if and only if $\rho$ lifts to $\GL_n(R)$.
	
	In particular, let $p$ be a prime, let $\widetilde V$ be a free $(\Z/p^2\Z)$-module of finite rank, and set $V\coloneqq\widetilde V/p\widetilde V$. Then \eqref{eq:square-zero} becomes
	\begin{equation}\label{eq:lifting-extension}
		1\longrightarrow\End_{\F_p}(V)
		\xlongrightarrow{\iota_{\widetilde V}}\GL(\widetilde V)
		\longrightarrow\GL(V)\longrightarrow1,
	\end{equation}
	where $\iota_{\widetilde V}(A)\coloneqq1_{\widetilde V}+p\widetilde A$ for any lift $\widetilde A\in\End_{\Z/p^2\Z}(\widetilde V)$ of $A\in\End_{\F_p}(V)$. The map $\iota_{\widetilde V}$ does not depend on the choice of $\widetilde A$, and the induced action on $\End_{\F_p}(V)$ is the conjugation action. For every continuous homomorphism $\rho\colon G\to\GL(V)$, we also write $o(\rho)\coloneqq o_{\Z/p^2\Z}(\rho)\in H^2(G,\ad(\rho))$ for the pullback class of \eqref{eq:lifting-extension}.

	\subsection{The Witt-vector extension}\label{subsec:witt}
	
	Let $k$ be a field of characteristic $p>0$, and let $W_2(k)$ be the ring of $p$-typical Witt vectors of length $2$. As a set, $W_2(k)=k\times k$, and the Witt-vector operations satisfy
	\begin{equation}\label{eq:witt-identities}
		(0,x)+(0,y)=(0,x+y),\quad (0,x)(0,y)=0,\quad (a,b)(0,x)=(0,a^p x).
	\end{equation}
	In particular, we have an exact sequence of additive groups
	\begin{equation}\label{eq:k-w-k}
		0\longrightarrow k\longrightarrow W_2(k)\longrightarrow k\longrightarrow0,
	\end{equation}
	where the first map sends $b$ to $(0,b)$ and the reduction map sends $(a,b)\mapsto a$; see, for example, \cite[(2.3)]{MSMA26}.
	
	Let $n\geq1$, and write $\mathfrak{gl}_n$ for the Lie algebra of $\mathrm{GL}_n$, so that we have an identification $\mathfrak{gl}_n(k)=M_n(k)$ with the conjugation action of $\GL_n(k)$. We denote by $\mathfrak{gl}_n^{(1)}$ the Frobenius twist of this representation, so that $g\in\GL_n$ acts by $y\mapsto \mathrm{Fr}(g)y\mathrm{Fr}(g)^{-1}$, where $\mathrm{Fr}$ raises the matrix entries to the $p$th power.
	
	For every commutative $k$-algebra $A$, let $\widetilde{\GL}_n(A)\coloneqq\GL_n(W_2(A))$. Since $W_2(A)=A\times A$ as sets and the Witt-vector operations are given by polynomials, the functor $\widetilde{\GL}_n$ is represented by an affine group scheme over $k$, whose underlying scheme is $\GL_n\times\mathbb A^{n^2}_{k}$. Moreover, \eqref{eq:k-w-k} induces a short exact sequence
	\begin{equation}\label{eq:witt-extension}
		0\longrightarrow\mathfrak{gl}_n^{(1)}
		\longrightarrow\widetilde{\GL}_n
		\xlongrightarrow{\pi}\GL_n\longrightarrow1,
	\end{equation}
	where the first map sends $y=(y_{ij})$ to $I+((0,y_{ij}))$. We denote by 
	\begin{equation}\label{eq:universal-obstruction}
		e_n\in H^2_{\mathrm{rat}}(\GL_n,\mathfrak{gl}_n^{(1)})
	\end{equation}
	the class of this extension. Taking $k$-points in \eqref{eq:witt-extension}, we obtain the exact sequence
	\begin{equation}\label{eq:witt-square-zero}
		0\longrightarrow M_n(k)^{(1)}
		\longrightarrow\GL_n(W_2(k))
		\longrightarrow\GL_n(k)\longrightarrow1.
	\end{equation}
	
	Let $G$ be an abstract or profinite group, let $V$ be a finite-dimensional $G$-representation over $k$, and let $\rho\colon G\to\GL_n(k)$ be the corresponding homomorphism after choosing a basis of $V$. By \eqref{eq:witt-square-zero}, reduction modulo $p$ defines an obstruction $o(V)\in H^2(G,\End(V)^{(1)})$, which vanishes if and only if $V$ lifts to a $W_2(k)$-free $G$-module. By functoriality, we have
	\begin{equation}\label{eq:obstruction-pullback}
		o(V)=\rho^*e_n\quad \text{in $H^2(G,\End(V)^{(1)})$}.
	\end{equation}

	\subsection{Reduction to the Witt-vector extension}
	
	\begin{lemma}\label{lem:reduction}
		Let $G$ be an abstract or profinite group, let $p$ be a prime, and let $n\geq1$. Suppose that either
		\begin{enumerate}
			\item for every algebraically closed field $k$ of characteristic $p$, every continuous homomorphism $G\to\GL_n(k)$ lifts to $\GL_n(W_2(k))$; or
			\item $G$ is profinite, $H^1(U,\F_p)$ is finite for every open normal subgroup $U$ of $G$, and every continuous homomorphism $G\to\GL_n(\overline\F_p)$ lifts to $\GL_n(W_2(\overline\F_p))$.
		\end{enumerate}
		Then $G$ is liftable in characteristic $p$.
	\end{lemma}

	\begin{proof}
		We first show that if a profinite group $G$ satisfies the assumptions of (2), then it satisfies the assumptions of (1). Let $k$ be an algebraically closed field of characteristic $p$, let $\rho\colon G\to\GL_n(k)$ be a continuous homomorphism, and set $U\coloneqq\ker(\rho)$ and $U_0\coloneqq\bigcap_{\chi\in H^1(U,\F_p)}\ker(\chi)$. By assumption (2), $U_0$ is open and characteristic in $U$, hence normal in $G$.
		
		Since the groups $G/U$ and $G/U_0$ are finite, we may consider the affine $\overline\F_p$-scheme of finite type $X\coloneqq\operatorname{Hom}(G/U,\GL_n)$, and let $\widetilde X$ be the closed subscheme of the affine $\overline \F_p$-scheme of finite type $\operatorname{Hom}(G/U_0,\widetilde{\GL}_n)$ parametrizing homomorphisms whose reductions factor through $G/U$. The reduction map $\widetilde{\GL}_n\to \GL_n$ induces a morphism $f\colon\widetilde X\to X$, and we claim that $f$ is surjective. For this, it suffices to show that $f$ is surjective on $\overline \F_p$-points. Every $\overline\F_p$-point of $X$ corresponds to a homomorphism $G/U\to\GL_n(\overline\F_p)$, which lifts by assumption (2) to a continuous homomorphism $G\to\GL_n(W_2(\overline\F_p))$. Since the kernel of $\GL_n(W_2(\overline\F_p))\to\GL_n(\overline\F_p)$ is an elementary abelian $p$-group, this lift sends $U_0$ to the identity and hence factors through $G/U_0$. Thus $f$ is surjective on $\overline \F_p$-points and hence it is surjective, as claimed. 
		
		Now fix an embedding $\overline\F_p\hookrightarrow k$, and observe that the fiber over the $k$-point of $X$ defined by $\rho$ is a nonempty scheme of finite type over the algebraically closed field $k$, and so has a $k$-point, which then gives a lift of $\rho$ factoring through $G/U_0$, hence in particular a continuous lift of $\rho$ proving (1).
		
		Now suppose that (1) holds. Let $R$ be a commutative local ring with maximal ideal $\m$ such that $\m^2=0$, let $k\coloneqq R/\m$ have characteristic $p$, and let $\rho\colon G\to\GL_n(k)$ be a continuous homomorphism. By \cite[Lemma~03C3]{Stacks}, there is a flat local homomorphism $R\to R'$ such that $\m R'$ is the maximal ideal of $R'$ and $R'/\m R'$ is an algebraic closure $\bar{k}$ of $k$. Since $R\to R'$ is flat, we have $\m R'=\m\otimes_k\bar{k}$, and hence $(\m R')^2=0$. The map
		\[
		H^2(G,M_n(\m))\longrightarrow
		H^2(G,M_n(\m)\otimes_k\bar{k})
		=H^2(G,M_n(\m R'))
		\]
		sends $o_R(\rho)$ to $o_{R'}(\rho\otimes_k\bar{k})$ and is split injective: a $k$-linear retraction of the inclusion $k\hookrightarrow\bar{k}$ induces a splitting. It therefore suffices to show that $o_{R'}(\rho\otimes_k\bar{k})=0$.
		
		Since $\bar{k}$ is perfect and $(\m R')^2=0$, there is a ring map $W_2(\bar{k})\to R'$ inducing the identity on residue fields; see \cite[Chapter~II, \S5, Proposition~10]{Ser79}. By (1), the representation $\rho\otimes_k\bar{k}$ lifts to $\GL_n(W_2(\bar{k}))$, and hence to $\GL_n(R')$. Thus $o_{R'}(\rho\otimes_k\bar{k})=0$, as desired.
	\end{proof}

	\subsection{Further reduction lemmas}
	
	The following observation was suggested by Mathieu Florence.
	
	\begin{lemma}\label{lem:two-coefficients}
		Let $p$ be a prime, let $k$ be a finite extension of $\F_p$, let $G$ be a profinite group, and let $\rho\colon G\to\operatorname{GL}_n(k)$ be a continuous homomorphism.
		Then $\rho$ lifts to $\operatorname{GL}_n(W_2(k))$ if and only if its restriction of scalars from $k$ to $\F_p$ lifts to $\operatorname{GL}_{n[k:\F_p]}(\Z/p^2\Z)$.
	\end{lemma}
	
	\begin{proof}
		Set $r\coloneqq[k:\F_p]$, let $V$ be the $k$-representation defined by $\rho$, and let $V_0$ be the restriction of scalars of $V$ to $\F_p$.
		Since $W_2(k)$ is free of rank $r$ over $\Z/p^2\Z$, a lift of $V$ to $W_2(k)$ gives a lift to $\Z/p^2\Z$ of $V_0$.
		
		Conversely, suppose that $V_0$ lifts to a free $\Z/p^2\Z$-module $\widetilde V_0$.
		Then the free $W_2(k)$-module $W_2(k)\otimes_{\Z/p^2\Z}\widetilde V_0$ lifts the $k$-representation $k\otimes_{\F_p}V_0$.
		Since $k/\F_p$ is finite and Galois, we have an isomorphism of $k$-algebras
		\[
		k\otimes_{\F_p}k\simeq\prod_{\sigma\in\operatorname{Gal}(k/\F_p)}k,
		\qquad a\otimes b\longmapsto(a\sigma(b))_\sigma.
		\]
		The factor corresponding to the identity automorphism shows that $V$ is a direct summand of $k\otimes_{\F_p}V_0$ as a $k$-representation of $G$.
		By \cite[Lemma~3.4]{DCF22}, the direct summand $V$ lifts to $W_2(k)$, as desired.
	\end{proof}
	
	The following observation will allow us to reduce to the case when $\mu_p\subset F$ in the proof of Theorem~\ref{thm:main-odd}.
	
	\begin{lemma}\label{lem:induction}
		Let $p$ be a prime, let $\Lambda$ be either $\Z_p$ or $\Z/p^e\Z$ for some integer $e\geq1$, let $G$ be a profinite group, let $H\subset G$ be an open subgroup of index $m$, and let $\rho'\colon H\to\operatorname{GL}_n(\F_p)$ be a continuous homomorphism.
		Set $\rho\coloneqq\operatorname{Ind}_H^G\rho'$.
		\begin{enumerate}[label=\textup{(\arabic*)}]
			\item\label{ind:first} If $\rho'$ lifts to $\operatorname{GL}_n(\Lambda)$, then $\rho$ lifts to $\operatorname{GL}_{mn}(\Lambda)$.
			\item\label{ind:second} Let $U\subset G$ be a closed subgroup.
			If the restriction of $\rho'$ to $H\cap g^{-1}Ug$ lifts to $\operatorname{GL}_n(\Lambda)$ for every $g\in G$, then $\rho|_U$ lifts to $\operatorname{GL}_{mn}(\Lambda)$.
			\item\label{ind:third} If $\rho$ lifts to $\operatorname{GL}_{mn}(\Z/p^2\Z)$, then $\rho'$ lifts to $\operatorname{GL}_n(\Z/p^2\Z)$.
		\end{enumerate}
	\end{lemma}
	
	\begin{proof}
		Part (1) follows from the compatibility of induction with reduction modulo $p$.
		Part (2) follows from (1) and Mackey's formula \cite[\S7.3, Proposition~22]{Ser77}.
		For (3), the restriction $\rho|_H$ contains $\rho'$ as a direct summand.
		Since a continuous lift to $\Z/p^2\Z$ factors through a finite quotient, the conclusion follows from \cite[Lemma~2.2]{KM26}.
	\end{proof}

	\section{Quadratic obstruction classes in characteristic two}\label{sec:quadratic}

	\subsection{An exact sequence}\label{subsec:quadratic}
	
	Let $k$ be a field of characteristic $2$, let $\mathcal G$ be an affine group scheme over $k$, and let $W$ be a finite-dimensional representation of $\mathcal G$. We have an exact sequence of $\mathcal G$-representations
	\begin{equation}\label{eq:quadratic-sequences}
		0 \longrightarrow W^{(1)}\xrightarrow{\ \iota \ } \Sym^2(W)
		\xrightarrow{\ N\ }W\otimes W
		\xrightarrow{\ \pi\ }\Sym^2(W)\longrightarrow0,
	\end{equation}
	where, for all $w,w'\in W$, the map $\iota$ sends $\mathrm{Fr}(w)$ to $w^2$, the map $N$ sends $ww'$ to $w\otimes w'+w'\otimes w$, and $\pi$ sends $w\otimes w'$ to $ww'$.
	
	For every $q\in\Sym^2(W)^{\mathcal G}$, the double connecting homomorphism associated with \eqref{eq:quadratic-sequences} sends $q$ to a class
	\[
		\kappa_q\in H^2_{\mathrm{rat}}(\mathcal G,W^{(1)}).
	\]
	Let $G$ be an abstract or profinite group, and let $W$ be a finite-dimensional representation of $G$ over $k$. For every $q\in(\Sym^2W)^G$, the same sequence defines a class 
	\[\kappa_q\in H^2(G,W^{(1)}),\] whose formation is compatible with pushforward along homomorphisms of $G$-represen\-tations $W\to W'$ over $k$, pullback along homomorphisms $G'\to G$, and extension of scalars.
	
	The class $\kappa_q$ may be refined as follows. For every exact sequence of $G$-representations
	\begin{equation}\label{eq:k-e-w}
	0\longrightarrow k\longrightarrow E\xlongrightarrow{\eta}W\longrightarrow0,
	\end{equation}
	we have an exact sequence
	\begin{equation}\label{eq:quadratic-lifted}
		0\longrightarrow E^{(1)}\xlongrightarrow{\widetilde \iota}\Sym^2E
		\xlongrightarrow{\widetilde N}W\otimes E
		\xlongrightarrow{\widetilde\pi}\Sym^2W\longrightarrow0,
	\end{equation}
	where, for all $e,e'\in E$ and $w\in W$, the map $\widetilde\iota$ sends $\mathrm{Fr}(e)$
	to $e^2$, the map $\widetilde N$ sends $ee'$ to $\eta(e)\otimes e'+\eta(e')\otimes e$, and $\widetilde\pi$ sends $w\otimes e$ to $w\eta(e)$. To prove that \eqref{eq:quadratic-lifted} is exact, it suffices to write $E=kz\oplus W$, where $z$ generates the kernel of $p$, and then to apply the exactness of \eqref{eq:quadratic-sequences}. For every $q\in(\Sym^2W)^G$, the double connecting homomorphism for \eqref{eq:quadratic-lifted} sends $q$ to
	\[
	\widetilde\kappa_{q,E}\in H^2(G,E^{(1)}).
\]
	
	\begin{lemma}\label{lem:lifting-kappa}
		Let $G$ be an abstract or profinite group, let $k$ be a field of characteristic $2$, let $W$ be a finite-dimensional $G$-representation over $k$. For every exact sequence of the form \eqref{eq:k-e-w} and every $q\in(\Sym^2W)^G$, we have
		\[
		 \eta^{(1)}_*\widetilde\kappa_{q,E}=\kappa_q \quad\text{in }H^2(H,W^{(1)}).
		\]
		Moreover, for every homomorphism of abstract groups or of profinite groups $\psi\colon H\to G$ and for every $H$-equivariant section $s\colon W\to E$ of $\eta$, we have
		\[
		\psi^*\widetilde\kappa_{q,E}=s^{(1)}_*\psi^*\kappa_q
		\quad\text{in }H^2(H,E^{(1)}).
		\]
	\end{lemma}
	\begin{proof}
		The equality $p^{(1)}_*\widetilde\kappa_{q,E}=\kappa_q$ follows from the following commutative diagram with exact rows
		\[
		\begin{tikzcd}
			0 \arrow[r] &
			E^{(1)} \arrow[r, "\widetilde \iota"] \arrow[d, "p^{(1)}"] &
			\Sym^2 E \arrow[r, "\widetilde N"] \arrow[d, "\Sym^2 p"] &
			W\otimes E \arrow[r, "\widetilde\pi"] \arrow[d, "1_W\otimes p"] &
			\Sym^2 W \arrow[r] \arrow[d, equal] &
			0 \\
			0 \arrow[r] &
			W^{(1)} \arrow[r, "\iota"] &
			\Sym^2 W \arrow[r, "N"] &
			W\otimes W \arrow[r, "\pi"] &
			\Sym^2 W \arrow[r] &
			0
		\end{tikzcd}
		\]
		and the naturality of connecting homomorphisms.
		
		After pullback along $\psi$, the maps $s^{(1)}$, $\Sym^2s$, $1_W\otimes s$ and the identity of $\Sym^2W$ give a morphism from \eqref{eq:quadratic-sequences} to \eqref{eq:quadratic-lifted}, and the conclusion again follows from the naturality of connecting maps.
	\end{proof}

	\begin{prop}\label{prop:annihilation}
		Let $G$ be an abstract or profinite group, let $k$ be a field of characteristic $2$, and let $W$ be a finite-dimensional $G$-representation over $k$. For every $q\in(\Sym^2 W)^G$ and every $\xi\in H^1(G,(W^{(1)})^\vee)$, we have
		\[
		\langle\kappa_q\cup\xi\rangle=0\quad\text{in }H^3(G,k),
		\]
		where $\langle-,-\rangle$ denotes the evaluation pairing $W^{(1)}\times(W^{(1)})^\vee\to k$.
	\end{prop}
	
	\begin{proof}
		Since the formation of $\kappa_q$ and cup products commute with extension of scalars and $H^3(G,k)\to H^3(G,\bar k)$ is split injective, we may assume that $k$ is algebraically closed, and hence in particular perfect.
		Fix $\xi\in H^1(G,(W^{(1)})^\vee)$. Via the canonical identification
		\[
		H^1(G,(W^{(1)})^\vee)\simeq\operatorname{Ext}^1_G(W^{(1)},k),
		\]
		the class $\xi$ corresponds to an extension of $W^{(1)}$ by $k$. Since $k$ is perfect, Frobenius twisting is an exact auto-equivalence of the category of $G$-representations over $k$, so this extension is of the form
		\[
		0\longrightarrow k\longrightarrow E^{(1)}\xlongrightarrow{\eta^{(1)}}W^{(1)}\longrightarrow0,
		\]
		for some extension of the form \eqref{eq:k-e-w}. We let
		\[
		\partial_\xi\colon H^2(G,W^{(1)})\longrightarrow H^3(G,k)
		\]
		be the connecting homomorphism associated with this extension, which is given by cup product with the extension class $\xi$:
		\[
		\partial_\xi(c)=\langle c\cup\xi\rangle\qquad\text{for every }c\in H^2(G,W^{(1)}).
		\]
		By Lemma~\ref{lem:lifting-kappa}, the class $\kappa_q$ lifts to $H^2(G,E^{(1)})$, and hence the cohomology long exact sequence gives
		\[
		\langle\kappa_q\cup\xi\rangle=\partial_\xi(\kappa_q)=0.\qedhere
		\]
	\end{proof}

	\subsection{\texorpdfstring{Comparison of the two classes $e_n$ and $\kappa_n$}{Comparison of the two classes e\_n and kappa\_n}}\label{subsec:comparison}

	Let $n\ge2$. For $x\in\mathfrak{gl}_n$, write
	\[
	\det(TI_n-x)
	=T^n-c_1(x)T^{n-1}+c_2(x)T^{n-2}
	-\cdots+(-1)^n c_n(x),
	\]
	for some $c_i\in k[\mathfrak{gl}_n]_i=\Sym^i(\mathfrak{gl}_n^\vee)$. Since the characteristic polynomial is invariant under conjugation, we have $c_i\in \Sym^i(\mathfrak{gl}_n^\vee)^{\mathrm{GL}_n}$.
	Applying the construction of
	Subsection~\ref{subsec:quadratic} to $\mathcal{G}=\mathrm{GL}_n$, $W=\mathfrak{gl}_n^\vee$ and $q=c_2$, we obtain a class $\kappa_{c_2}\in H^2_{\mathrm{rat}}
	(\GL_n,(\mathfrak{gl}_n^{\vee})^{(1)})$. 
	Via the isomorphism $(\mathfrak{gl}_n^\vee)^{(1)}\xrightarrow{\sim}\mathfrak{gl}_n^{(1)}$
	induced by the trace pairing, we regard $\kappa_{c_2}$ as a class
	\[\kappa_n\in H^2_{\mathrm{rat}}(\GL_n,\mathfrak{gl}_n^{(1)}).\]
	For $n=1$, we set $c_2=0$ and $\kappa_1=0$.

	\begin{prop}\label{prop:comparison}
		Let $k$ be a field of characteristic $2$ and let $n\ge1$. Then
		\[
		e_n=\kappa_n
		\quad\text{in}\quad
		H^2_{\mathrm{rat}}(\GL_n,\mathfrak{gl}_n^{(1)}).
		\]
	\end{prop}

	\begin{proof}
		For $n=1$, the Teichm\"uller section $\GL_1\to\widetilde{\GL}_1$ shows that $e_1=0$ and $\kappa_1=0$. From now on, we assume $n\ge2$. Since both classes are defined over $\mathbb F_2$ and rational cohomology commutes with extension of the ground field \cite[Proposition~I.4.18]{Jan03}, we may also assume that $k=\mathbb F_2$.
		
		By \cite[Remark~5.2]{TVK10},
		\[
		\dim_k H^2_{\mathrm{rat}}(\GL_n,\mathfrak{gl}_n^{(1)})=1,\qquad e_n\neq 0.
		\]
		It thus suffices to prove that $\kappa_n\neq 0$. For this, consider the embedding $\GL_2\to\GL_n$ given by $g\mapsto\operatorname{diag}(g,I_{n-2})$. The inclusion $\mathfrak{gl}_2\hookrightarrow\mathfrak{gl}_n$ into the upper-left block induces a $\GL_2$-equivariant map $\mathfrak{gl}_n^\vee\to\mathfrak{gl}_2^\vee$ which sends $c_2$ to the determinant upon taking second symmetric powers. By naturality of the exact sequence \eqref{eq:quadratic-sequences}, it suffices to prove that $\kappa_2\neq 0$.
		
		Consider the $k$-group embedding
		\[
		\alpha_2\xhookrightarrow{\phantom{aaa}}\GL_2,\qquad
		\lambda\longmapsto
		\begin{pmatrix}1&\lambda\\0&1\end{pmatrix}.
		\]
		Since the distribution algebra of $\alpha_2$ is $R=k[D]/(D^2)$, for every rational $\alpha_2$-module $M$ we have $H^i_{\mathrm{rat}}(\alpha_2,M)\simeq\operatorname{Ext}^i_R(k,M)$; 
		see \cite[I.4 and I.8]{Jan03}. The trivial $R$-module $k$ has the free resolution
		\[
		\cdots\xlongrightarrow{D}R\xlongrightarrow{D}R\xlongrightarrow{D}R\longrightarrow k\longrightarrow0.
		\]
		Applying $\Hom_R(-,M)$, for any $R$-module $M$, gives the complex
		\begin{equation}\label{eq:alpha-periodic}
			M\xlongrightarrow{D}M\xlongrightarrow{D}M\xlongrightarrow{D}\cdots
		\end{equation}
		Since $\alpha_2$ acts trivially on $(\mathfrak{gl}_2^\vee)^{(1)}$, this gives an isomorphism
		\begin{equation}\label{eq:h2rat-w1}
			H^2_{\mathrm{rat}}(\alpha_2,(\mathfrak{gl}_2^\vee)^{(1)})\simeq(\mathfrak{gl}_2^\vee)^{(1)}.
		\end{equation}
		Let $a,b,c,d\in\mathfrak{gl}_2^\vee$ be the coordinate functions defined by
		\[
		x=\begin{pmatrix}a(x)&b(x)\\c(x)&d(x)\end{pmatrix}
		\]
		for every $k$-algebra $A$ and every $x\in \mathfrak{gl}_2(A)$. Then $c_2=ad+bc$, and the action of $D$ on $\mathfrak{gl}_2^\vee$ is given by
		\[
		D(a)=D(d)=c,\qquad D(b)=a+d,\qquad D(c)=0.
		\]
		Define
		\[
		B\coloneqq a\otimes d+b\otimes c\in\mathfrak{gl}_2^\vee\otimes\mathfrak{gl}_2^\vee,
		\]
		and observe that $\pi(B)=ad+bc=c_2$ and that the Leibniz rule gives
		\[
		D(B)=c\otimes d+d\otimes c=N(cd),
		\qquad D(cd)=c^2=\iota(\mathrm{Fr}(c)).
		\]
		We compute the two successive connecting homomorphisms using \eqref{eq:alpha-periodic}: using that $N(cd)=D(B)$ and $\iota(\mathrm{Fr}(c))=D(cd)$, the first connecting map sends $c_2$ to the class of $D(B)$, while the second sends this class to $\mathrm{Fr}(c)$. All in all, this shows that the restriction of $\kappa_2$ to $\alpha_2$ corresponds to $\mathrm{Fr}(c)\neq 0$ under \eqref{eq:h2rat-w1}, so that in particular $\kappa_2\neq 0$, as desired.
	\end{proof}

	\begin{proof}[Proof of Theorem~\ref{thm:universal}]
		Let $n\coloneqq\dim_kV$, and let $\rho\colon G\to\GL_n(k)$ be a homomorphism corresponding to $V$. By \eqref{eq:obstruction-pullback} and Proposition~\ref{prop:comparison}, we have $o(V)=\rho^*\kappa_n$. Therefore, applying Proposition~\ref{prop:annihilation} to $W=\End(V)^\vee$ and $q=c_2$, we obtain $\tr_*(o(V)\cup\xi)=0$ for every $\xi\in H^1(G,\End(V)^{(1)})$, as desired.
	\end{proof}

	Theorem~\ref{thm:universal} has the following consequence for the lifting problem.
	
	\begin{cor}\label{cor:perfect-implies-liftable}
		Let $G$ be an abstract or profinite group, and assume that, for every algebraically closed field $k$ of characteristic $2$ and for every finite-dimensional $G$-representation $W$ over $k$ equipped with a nondegenerate $G$-invariant bilinear pairing $\langle-,-\rangle\colon W\times W\to k$, the pairing \[H^2(G,W)\times H^1(G,W)\longrightarrow H^3(G,k),\qquad (x,y)\longmapsto\langle x\cup y\rangle,\]
		is nondegenerate in the first variable, that is, the induced map
		\[H^2(G,W)\longrightarrow\operatorname{Hom}_k(H^1(G,W),H^3(G,k))\]
		is injective. Then $G$ is liftable in characteristic $2$.
	\end{cor}
	
	\begin{proof}
		Let $k$ be an algebraically closed field of characteristic $2$, and let $V$ be a finite-dimensional 	$G$-representation over $k$. By Theorem~\ref{thm:universal}
		applied to the trace pairing on $\End(V)^{(1)}$, which is nondegenerate and $G$-invariant, we have $o(V)=0$ in $H^2(G,\End(V)^{(1)})$, so that $V$ lifts to $W_2(k)$. The conclusion follows from Lemma~\ref{lem:reduction}.
	\end{proof}
	
	\begin{rmk}
		In Corollary~\ref{cor:perfect-implies-liftable}, it suffices to assume nondegeneracy in the first variable for the cup-product pairing
		\[H^2(G,\End(V)^{(1)})\times H^1(G,\End(V)^{(1)})\longrightarrow H^3(G,k)\]
		induced by the trace pairing, for every algebraically closed field $k$ of characteristic $2$ and every finite-dimensional $G$-representation $V$ over $k$.
	\end{rmk}
	
	\section{Proof of Theorem~\ref{thm:quaternion}}\label{sec:ordinary-duality}

	\subsection{Generalized quaternion groups}
	The following lemma is certainly well known, but we include a proof since we could not find a precise reference. 
	
	\begin{lemma}\label{lem:duality}
		Let $m\geq2$ be a power of $2$, let $k$ be a field of characteristic $2$, and let $N$ be a finite-dimensional $k[Q_{4m}]$-module equipped with a nondegenerate $Q_{4m}$-invariant bilinear pairing $\langle-,-\rangle\colon N\times N\to k$. Then $H^3(Q_{4m},k)\simeq k$, and the pairing
		\[
		H^2(Q_{4m},N)\times H^1(Q_{4m},N)\longrightarrow H^3(Q_{4m},k)\simeq k,\qquad (x,y)\longmapsto\langle x\cup y\rangle,
		\]
		is perfect.
	\end{lemma}
	
	\begin{proof}
		The group $Q_{4m}$ is isomorphic to a subgroup of the group of unit quaternions, and hence acts freely on the $3$-sphere $S^3$. It follows from \cite[Chapter~VI, \S9]{Brown82} that the cohomology of $Q_{4m}$ is $4$-periodic, that is, there exists $u\in\hat H^4(Q_{4m},\Z)$ such that, for every $k[Q_{4m}]$-module $N$ and every $i\in\Z$, cup product with $u$ induces isomorphisms of Tate cohomology groups
		\[
		u\cup(-)\colon\hat H^i(Q_{4m},N)\xlongrightarrow{\sim}\hat H^{i+4}(Q_{4m},N).
		\]
		By Tate duality \cite[Chapter~VI, \S7]{Brown82} and the identification $N\simeq N^\vee$ induced by $\langle-,-\rangle$, the pairing
		\[
		H^2(Q_{4m},N)\times\hat H^{-3}(Q_{4m},N)\longrightarrow\hat H^{-1}(Q_{4m},k)
		\]
		induced by cup product and $\langle-,-\rangle$ is perfect. Cup product with $u$ gives isomorphisms
		\[
		\hat H^{-3}(Q_{4m},N)\xlongrightarrow{\sim}H^1(Q_{4m},N),\qquad
		\hat H^{-1}(Q_{4m},k)\xlongrightarrow{\sim}H^3(Q_{4m},k),
		\]
		fitting in the following commutative square
		\[
		\begin{tikzcd}[column sep=large, row sep=large]
			H^2(G,N)\times\widehat H^{-3}(G,N)
			\arrow[r, "\cup"]
			\arrow[d, "{\mathrm{id}\times(u\cup-)}", "\wr"']
			&
			\widehat H^{-1}(G,k)
			\arrow[d, "u\cup-", "\wr"']
			\\
			H^2(G,N)\times H^1(G,N)
			\arrow[r, "\cup"]
			&
			H^3(G,k).
		\end{tikzcd}
		\]
		Hence the pairing in the statement is perfect. Finally, since $\operatorname{char}(k)=2$ and $|Q_{4m}|=4m$, by the $4$-periodicity property we have
		\[
		H^3(Q_{4m},k)\simeq\hat H^{-1}(Q_{4m},k)=k.\qedhere
		\]
	\end{proof}

	\begin{rmk}
		Lemma~\ref{lem:duality} can also be deduced from Poincar\'e duality with local coefficients for the closed oriented $3$-manifold $S^3/Q_{4m}$; see \cite[\S3.H]{Hat02}. Indeed, letting $\mathcal N$ be the local system associated with $N$, the classifying map $S^3/Q_{4m}\to BQ_{4m}$ induces isomorphisms $H^i(Q_{4m},N)\xrightarrow{\sim}H^i(S^3/Q_{4m},\mathcal N)$ for $i=1,2$ and $H^3(Q_{4m},k)\xrightarrow{\sim}H^3(S^3/Q_{4m},k)$. The assertion then follows from Poincar\'e duality.
	\end{rmk}

	\begin{proof}[Proof of Theorem~\ref{thm:quaternion}]
		Let $R$ be a commutative local ring with maximal ideal $\m$ such that $\m^2=0$, and let $k\coloneqq R/\m$. If $\operatorname{char}(k)\neq 2$, then $|Q_{4m}|$ is invertible in $k$, and hence every representation of $Q_{4m}$ over $k$ lifts to $R$. Suppose that $\operatorname{char}(k)=2$. For every algebraically closed field $k'$ of characteristic $2$, Lemma~\ref{lem:duality} shows that $H^3(Q_{4m},k')\simeq k'$ and that the pairings of Corollary~\ref{cor:perfect-implies-liftable} are perfect, so that $Q_{4m}$ is liftable in characteristic $2$, as desired.
	\end{proof}

	\begin{proof}[Proof of Corollary~\ref{cor:classification}]
		If $G$ is liftable, then it is $p$-liftable for every prime $p$, and hence every subgroup of $G$ is $p$-liftable for every prime $p$. By \cite[Corollary~3.2]{KM26v1}, this implies that $G$ has no element of prime order $p\geq5$ and no element of order $9$. Moreover, by \cite[Proposition~5.1]{KM26v1} (resp. \cite[Proposition~4.1]{KM26v1}), $G$ has no subgroup isomorphic to $C_3\times C_3$ (resp. $C_2\times C_2$). Since every $3$-group of order at least $9$ contains a subgroup of order $9$, which is isomorphic to either $C_9$ or $C_3\times C_3$, the Sylow $3$-subgroups of $G$ have order at most $3$, proving (1). Let $P$ be a nontrivial Sylow $2$-subgroup of $G$, and choose a central involution $z\in P$. If $z'$ is an involution distinct from $z$, then $P$ contains $\langle z,z'\rangle\simeq C_2\times C_2$, a contradiction. Thus $P$ has a unique subgroup of order $2$, and hence is cyclic or generalized quaternion by \cite[Theorem~5.4.10(ii)]{Gor80}. This proves (2).
		
		Conversely, suppose that $G$ satisfies (1) and (2). Then every Sylow subgroup of $G$ is cyclic of $2$-power order, cyclic of order at most $3$, or generalized quaternion. Such groups are liftable by \cite[Corollary~3.4]{KM26v1} and Theorem~\ref{thm:quaternion}, and hence by \cite[Proposition~2.4]{KM26v1} so is $G$, as desired.
	\end{proof}

	\subsection{Two-dimensional local fields}
	Recall that a field $K$ is said to be a \emph{two-dimensional local field} if it is complete with respect to a discrete valuation whose residue field $K'$ is a nonarchimedean local field; if this is the case, the residue field of $K'$ is said to be the \emph{final residue field} of $K$. 
	
	\begin{cor}\label{cor:galois}
		Let $K$ be a two-dimensional local field whose final residue field has odd characteristic. Then the absolute Galois group of $K$ is liftable in characteristic $2$.
	\end{cor}
	
	For example, Corollary~\ref{cor:galois} applies to $K=\F_q(\!(u)\!)(\!(t)\!)$, where $q$ is odd, and to $K=F(\!(t)\!)$, where $F$ is a finite extension of $\Q_\ell$ for an odd prime $\ell$.

	\begin{proof}[Proof of Corollary~\ref{cor:galois}]
		Let $G_K$ be the absolute Galois group of $K$. By higher local duality \cite[Chapter~I, Theorem~2.17]{Mil06}, for every finite $2$-torsion $G_K$-module $A$, the cup product induces a perfect pairing
		\[
		H^2(G_K,A)\times H^1(G_K,\Hom(A,\F_2))
		\longrightarrow H^3(G_K,\F_2)\simeq\F_2,
		\]
		where we have used the fact that $\mu_2^{\otimes2}\simeq\F_2$ as a $G_K$-module. Using the trace pairing for finite extensions of $\F_2$, the same duality holds for finite-dimensional continuous representations over any finite field of characteristic $2$. Theorem~\ref{thm:universal} now implies that every continuous homomorphism $G_K\to\GL_n(\overline\F_2)$ lifts to $\GL_n(W_2(\overline\F_2))$.
		
		For every open subgroup $U$ of $G_K$, letting $L$ be the fixed field of $U$, by Kummer theory we have
		\[
		H^1(U,\F_2)=H^1(L,\F_2)\simeq L^\times/L^{\times2}\simeq(\Z/2\Z)^3.
		\]
		The conclusion follows from Lemma~\ref{lem:reduction}.
	\end{proof}

	\subsection{Three-manifolds}
	
	\begin{cor}\label{cor:manifold}
		Let $X$ be a closed connected $3$-manifold. Then the fundamental group $\pi_1(X)$ is liftable in characteristic $2$.
	\end{cor}

\begin{proof}
	Let $G\coloneqq\pi_1(X)$, let $k$ be an algebraically closed field of characteristic $2$, let $W$ be a finite-dimensional $G$-representation over $k$ equipped with a nondegenerate $G$-invariant bilinear pairing, and let $BG$ be the classifying space of $G$. Letting $\mathcal W$ be the local system on $X$ associated with $W$, the classifying map $c\colon X\to BG$ induces an isomorphism $H^1(G,W)\xrightarrow{\sim}H^1(X,\mathcal W)$ and an injection $H^2(G,W)\hookrightarrow H^2(X,\mathcal W)$.
	
	By Poincar\'{e} duality with local coefficients \cite[\S3.H]{Hat02}, the pairing
	\[H^2(X,\mathcal W)\times H^1(X,\mathcal W)\longrightarrow H^3(X,k)\simeq k\]
	is perfect. (We do not need to assume that the manifold $X$ is orientable because $\operatorname{char}(k)=2$.) It follows that the cup-product pairing
	\[H^2(G,W)\times H^1(G,W)\longrightarrow H^3(G,k)\]
	is nondegenerate in the first variable. The conclusion now follows from Corollary~\ref{cor:perfect-implies-liftable}.
\end{proof}
	
	\begin{rmk}
		Corollary~\ref{cor:manifold} is an analogue, in dimension $3$, of a result of Conti--Demarche--Florence \cite[Corollary~7.25]{CDF26}, which shows that, for every prime $p$, every $\F_p$-linear representation of the fundamental group of a closed connected orientable surface lifts to $\Z/p^2\Z$.
	\end{rmk}

	\section{Proof of Theorem~\ref{thm:main-even}}\label{sec:even-short}

	\subsection{Poitou--Tate duality}\label{sec:pt}
	
	Let $p$ be a prime, let $F$ be a global field of characteristic different from $p$ containing a primitive $p$th root of unity $\zeta$, let $S$ be a finite nonempty set of places of $F$ containing all archimedean places and all places above $p$, and let $G_{F,S}$ denote the Galois group of the maximal extension of $F$ unramified outside $S$.
	For every finite $p$-torsion $G_{F,S}$-module $N$ and every integer $i\geq0$, we set
	\[
	\Sh^i_S(F,N)\coloneqq
	\ker(H^i(G_{F,S},N)\longrightarrow\prod_{v\in S} H^i(F_v,N)).
	\]
	We let $N^D\coloneqq\operatorname{Hom}(N,\mu_p)$ be the Tate dual of $N$.
	By Poitou--Tate duality \cite[Chapter~I, Theorem~4.10(a)]{Mil06}, there is a perfect pairing
	\[
	\langle -,- \rangle_{\mathrm{PT}}\colon
	\Sh^1_S(F,N^D)\times\Sh^2_S(F,N)\longrightarrow\Q/\Z
	\]
	which can be explicitly described as follows; see \cite[Chapter~I, \S4, p.~65]{Mil06}.
	Let $x\in\Sh^1_S(F,N^D)$ and $y\in\Sh^2_S(F,N)$, and choose representative cocycles $f\in Z^1(G_{F,S},N^D)$ of $x$ and $\omega\in Z^2(G_{F,S},N)$ of $y$.
	For every place $v\in S$, choose $P_v\in N^D$ such that $\delta P_v=f_v$, and choose $h\in C^2(G_{F,S},\mu_p)$ such that $\delta h=f\cup\omega$. Such an $h$ exists because $[f\cup\omega]\in H^3(G_{F,S},\mu_p)$ is locally trivial at every place in $S$ and $\Sh^3_S(F,\mu_p)=0$;
	see \cite[Chapter~I, Theorem~4.10(c)]{Mil06}.
	Then for every place $v\in S$ we have $P_v\cup\omega_v-h_v\in Z^2(F_v,\mu_p)$, and
	\begin{equation}\label{eq:pt-definition}
		\langle x,y\rangle_{\mathrm{PT}}=\sum_{v\in S}\operatorname{inv}_v([P_v\cup\omega_v-h_v]).
	\end{equation}
	
\begin{lemma}\label{lem:pt-yoneda}
	Let $p$ be a prime, let $F$ be a global field of characteristic different from $p$ containing $\mu_p$, let $S$ be a finite nonempty set of places of $F$ containing all archimedean places and all places above $p$, let $M$ be a finite-dimensional $\F_p$-representation of $G_{F,S}$, and let $\xi\in\Sh^1_S(F,M^D)$ and $y\in\Sh^2_S(F,M)$. Let
	\begin{equation}\label{eq:pt-yoneda-extension}
		0\longrightarrow\mu_p\xlongrightarrow{j}N\xlongrightarrow{\pi}M\longrightarrow0
	\end{equation}
	be an extension of $\F_p[G_{F,S}]$-modules representing $\xi$. Suppose that $y$ admits a lift $\widetilde y\in H^2(G_{F,S},N)$, and for every place $v\in S$ choose a $G_{F_v}$-equivariant retraction $r_v\colon N\to\mu_p$ of $j$. Then
	\begin{equation}\label{eq:pt-yoneda-formula}
		\langle\xi,y\rangle_{\mathrm{PT}}
		=\sum_{v\in S}\operatorname{inv}_v((r_v)_*\widetilde y_v).
	\end{equation}
\end{lemma}

\begin{proof}
	We fix an $\F_p$-linear splitting of $\pi$ and use it to identify $N$ with $\mu_p\oplus M$ as an $\F_p$-vector space. Letting $f\in Z^1(G_{F,S},M^D)$ be a representative of $\xi$, the $G_{F,S}$-action on $N$ takes the form
	\[
	\sigma(c,m)=(c+\langle f(\sigma),\sigma m\rangle,\sigma m)
	\]
	for all $\sigma\in G_{F,S}$.
	
	We fix a cocycle $\omega\in Z^2(G_{F,S},M)$ representing $y$. Since $\pi_*\widetilde y=y$, we may choose $\widetilde\omega\in Z^2(G_{F,S},N)$ representing $\widetilde y$ such that $\pi\widetilde\omega=\omega$. Writing $\widetilde\omega=(-h,\omega)$ and letting $\mathrm{pr}_1\colon N\to\mu_p$ be the first projection, since $\widetilde\omega$ is a $2$-cocycle, for every $\sigma,\tau,\upsilon\in G_{F,S}$ we have
	\[
	\begin{aligned}
		0
		&=\operatorname{pr}_1(\delta\widetilde\omega
		(\sigma,\tau,\upsilon))\\
		&=-h(\tau,\upsilon)
		+\langle f(\sigma),\sigma\omega(\tau,\upsilon)\rangle
		+h(\sigma\tau,\upsilon)-h(\sigma,\tau\upsilon)+h(\sigma,\tau)\\
		&=-(\delta h)(\sigma,\tau,\upsilon)
		+(f\cup\omega)(\sigma,\tau,\upsilon),
	\end{aligned}
	\]
	that is, $\delta h=f\cup\omega$. Moreover, for every place $v\in S$, since $\xi_v=0$ we may choose $P_v\in M^D$ such that $f_v=\delta P_v$. This equality implies that the map
	$r_v\colon N\to\mu_p$ defined by $r_v(c,m)=c+\langle P_v,m\rangle$ is a $G_{F_v}$-equivariant retraction of $j$. Since $\widetilde y_v$ is represented by
	$\widetilde\omega_v=(-h_v,\omega_v)$, we have
	\[
	r_v\widetilde\omega_v=-h_v+\langle P_v,\omega_v\rangle=P_v\cup\omega_v-h_v\quad\text{in }Z^2(F_v,\mu_p).
	\]
	Now \eqref{eq:pt-yoneda-formula} follows from \eqref{eq:pt-definition} for the choice of $\widetilde y$ and $r_v$ made above.
	
	It remains to check that \eqref{eq:pt-yoneda-formula} is independent of the choices of $\widetilde y$ and $r_v$. For this, if $r_v'$ is another $G_{F_v}$-equivariant retraction of $j$, then $r_v'-r_v=\lambda_v \pi$ for some $G_{F_v}$-equivariant map $\lambda_v\colon M\to\mu_p$, and so since $y_v=0$ we have $(r_v')_*\widetilde y_v-(r_v)_*\widetilde y_v=(\lambda_v)_*y_v=0$. If $\widetilde y'$ is another lift of $y$, then $\widetilde y'-\widetilde y=j_*z$ for some $z\in H^2(G_{F,S},\mu_p)$, and since $r_vj=\mathrm{id}_{\mu_p}$ the two sums differ by $\sum_{v\in S}\operatorname{inv}_v(z_v)$, which vanishes by Brauer reciprocity and the fact that $z_v=0$ for every $v\notin S$.
\end{proof}

	\subsection{Proof of Theorem~\ref{thm:main-even}}
	The following proposition is the main input in the proof of Theorem~\ref{thm:main-even}. 
	
	\begin{prop}\label{prop:quadratic-global}
		Let $F$ be a global field of characteristic different from $2$, let $S$ be a finite nonempty set of places of $F$ containing all archimedean places and all places above $2$, and let $W$ be a finite-dimensional continuous $\F_2$-representation of $G_{F,S}$. Then, for every $q\in(\Sym^2 W)^{G_{F,S}}$, we have
		\[\kappa_q=0\quad \text{in $H^2(G_{F,S},W)$}.\]
	\end{prop}
	
	\begin{proof}
		Throughout this proof, we identify $\mu_2$ with $\F_2$ and hence $W^D$ with $W^\vee$. For a finite $2$-torsion $G_{F,S}$-module $N$, let
		\[
		\Sh^i_S(F,N)\coloneqq
		\ker (H^i(G_{F,S},N)\longrightarrow\prod_{v\in S}H^i(F_v,N)).
		\]
		We first show that $\kappa_q\in\Sh^2_S(F,W)$.  Let $v\in S$, and let $\lambda\in(W^D)^{G_{F_v}}$, which we regard as a $G_{F_v}$-equivariant map $\lambda\colon W\to\F_2$. By functoriality of the sequence \eqref{eq:quadratic-sequences} with respect to $\lambda$, we have $\lambda_*\kappa_{q,v}=\kappa_{(\Sym^2\lambda)(q)}$, where $\Sym^2\lambda\colon \Sym^2W\to \Sym^2\F_2$ is the map induced by $\lambda$. Since the map $\pi\colon\F_2\otimes\F_2\to\Sym^2\F_2$ in \eqref{eq:quadratic-sequences} is an isomorphism, the associated double connecting homomorphism is zero, and hence $\kappa_{(\Sym^2\lambda)(q)}=0$, that is, $\lambda_*\kappa_{q,v}=0$ for all $\lambda\in(W^D)^{G_{F_v}}$. 
		
		Suppose that $v$ is nonarchimedean. By local Tate duality \cite[Chapter~I, Corollary~2.3 and Theorem~2.13]{Mil06}, the pairing
		\[H^2(F_v,W)\times H^0(F_v,W^D)\longrightarrow H^2(F_v,\mu_2),\qquad (x,\lambda)\longmapsto\lambda_*x\]
		is perfect. Since $\lambda_*\kappa_{q,v}=0$ for every $\lambda\in H^0(F_v,W^D)$, we conclude that $\kappa_{q,v}=0$ if $v$ is nonarchimedean. If $v$ is real, the same argument
		applies replacing $H^0(F_v,W^D)$ by $\widehat H^0(F_v,W^D)$, which is a
		quotient of $(W^D)^{G_{F_v}}$, and if $v$ is complex then $H^2(F_v,W)=0$.
		Thus $\kappa_q\in\Sh^2_S(F,W)$, as desired. 
		
		Let $\xi\in\Sh^1_S(F,W^D)=\Sh^1_S(F,W^\vee)$, let
		\begin{equation}\label{eq:f2-e-w}
		0\longrightarrow\F_2\longrightarrow E\xlongrightarrow{\pi}W\longrightarrow0
		\end{equation}
		be an extension of $G_{F,S}$-representations over $\F_2$ corresponding to $\xi$, and let $\widetilde\kappa_{q,E}\in H^2(G_{F,S},B)$ be the canonical lift of $\kappa_q$ from Lemma~\ref{lem:lifting-kappa}. Since $\xi_v=0$, the extension has a $G_{F_v}$-equivariant splitting $s_v$ for every $v\in S$, and hence by Lemma~\ref{lem:lifting-kappa} applied to the homomorphism $G_{F_v}\to G_{F,S}$ we obtain $(\widetilde\kappa_{q,E})_v=(s_v)_*\kappa_{q,v}=0$, that is, $\widetilde\kappa_{q,E}\in\Sh^2_S(F,E)$.
		
		Dualizing \eqref{eq:f2-e-w} gives \[0 \longrightarrow W^D\xlongrightarrow{\pi^D}E^D \longrightarrow \F_2 \longrightarrow 0.\] Since $\xi$ is the image of $1\in \F_2= H^0(G_{F,S},\F_2)$ under the connecting map, the image of $\xi$ in $H^1(G_{F,S},E^D)$ is zero, and hence by the functoriality of the Poitou--Tate pairing
		\[
		\langle\xi,\kappa_q\rangle_{\mathrm{PT}}
		=\langle\xi,\pi_*\widetilde\kappa\rangle_{\mathrm{PT}}
		=\langle(\pi^D)_*\xi,\widetilde\kappa\rangle_{\mathrm{PT}}=0.
		\]
		Since the Poitou--Tate pairing $\Sh^1_S(F,W^D)\times\Sh^2_S(F,W)\to\Q/\Z$ is perfect, we obtain $\kappa_q=0$, as desired.
	\end{proof}

	\begin{thm}\label{thm:two-S}
		Let $F$ be a global field, let $S$ be a finite nonempty set of places of $F$, let $k$ be a field of characteristic $2$, and let $n$ be a positive integer. If $\operatorname{char}(F)\neq 2$, suppose that $S$ contains the archimedean places and all places above $2$. Then every continuous homomorphism $G_{F,S}\to\GL_n(k)$ lifts to $\GL_n(W_2(k))$.
	\end{thm}
	
	\begin{proof}
		Suppose first that $\operatorname{char}(F)=2$. Let $C$ be the smooth projective curve over the field of constants of $F$ with function field $F$, and set $U\coloneqq C\setminus S$. Since $S$ is nonempty, $U$ is affine, and $G_{F,S}$ is its \'{e}tale fundamental group. By \cite[Proposition~4.10]{DCF22}, we have $\operatorname{cd}_2(G_{F,S})\leq1$, that is, $H^2(G_{F,S},A)=0$ for every discrete $2$-primary torsion module $A$, so that in particular the lifting obstruction vanishes. From now on, we assume that $\operatorname{char}(F)\neq 2$.
		
		Suppose that $k=\F_2$, and let $V$ be the representation defined by a continuous homomorphism $\rho\colon G_{F,S}\to\GL_n(\F_2)$. By \eqref{eq:obstruction-pullback} and Proposition~\ref{prop:comparison}, the obstruction $o(V)$ is the class $\kappa_{c_2}$ associated to the second coefficient of the characteristic polynomial, with $W=\End(V)^\vee$ identified with $\End(V)$ by the trace pairing. By Proposition~\ref{prop:quadratic-global}, we have $o(V)=0$, and hence $\rho$ lifts to $\GL_n(\Z/4\Z)$. This proves Theorem~\ref{thm:two-S} when $k=\F_2$.
		
		By Lemma~\ref{lem:two-coefficients}, for every finite extension $\F_q/\F_2$, every continuous representation of $G_{F,S}$ over $\F_q$ lifts to $W_2(\F_q)$, and hence every continuous representation over $\overline\F_2$ lifts $W_2(\overline\F_2)$. For every open normal subgroup $U$ of $G_{F,S}$, there is a finite extension $L/F$ unramified outside $S$ such that $U=G_{L,S_L}$, where $S_L$ is the set of places of $L$ above $S$, and hence by \cite[Chapter~I, Corollary~4.15]{Mil06} the group $H^1(U,\F_2)$ is finite. Therefore Lemma~\ref{lem:reduction} applies, and we conclude that every continuous representation over $k$ lifts to $W_2(k)$, as desired.
	\end{proof}
	
	\begin{proof}[Proof of Theorem~\ref{thm:main-even}]
		Let $F$ be a global field, and let $\rho\colon G_F\to\GL_n(k)$ be a continuous homomorphism, and choose a finite nonempty set $S$ containing all places where $\rho$ is ramified and, if $\operatorname{char}(F)\neq 2$, all archimedean places and places above $2$. Then $\rho$ factors through $G_{F,S}$, and Theorem~\ref{thm:two-S} gives the desired lift.
	\end{proof}

	\section{Proof of Theorem~\ref{thm:main-odd}}\label{sec:odd-short}
	
	Throughout this section, we let $p$ be an odd prime.
	
\subsection{An observation on the trace map}

Let $G$ be a profinite group, let $k$ be a field, and let $V$ and $P$ be finite-dimensional $G$-representations over $k$. For every $G$-equivariant map $T\colon V\to P$, we write 
\[T_*\colon H^2(G,k)\longrightarrow H^2(G,\Hom_k(V,P))\] for the map induced by the $G$-equivariant homomorphism $k\to\Hom_k(V,P)$ sending $1$ to $T$.

The trace $\tr_V\colon\End_k(V)\to k$ is $G$-equivariant and hence induces a map
\[\Tr_V\colon H^2(G,\End_k(V))\longrightarrow H^2(G,k).\]
For every $T\in\End_G(V)$, the composite
\[H^2(G,k)\xlongrightarrow{T_*}H^2(G,\End_k(V))\xrightarrow{\Tr_V}H^2(G,k)\]
is given by multiplication by $\tr_V(T)$.
	
	\begin{lemma}\label{lem:criterion}
		Let $G$ be a profinite group, let $k$ be a field of characteristic $p>0$, and let $V$ be a finite-dimensional continuous $k$-representation of $G$. Let
		\begin{equation}\label{eq:detector}
			0\longrightarrow V\xlongrightarrow{i}P\xlongrightarrow{\pi}V\longrightarrow0
		\end{equation}
		be an extension, let $d\in\End_G(V)$, and let $\Phi\colon V\to P$ be equivariant with $\pi\Phi=d$. Let $B\in H^2(G,k)$, and suppose that $d_*B=0$ in $H^2(G,\End_k(V))$. Then there exists $\gamma\in H^2(G,k)$ such that, for every closed subgroup $H\subset G$ over which \eqref{eq:detector} splits, and every $H$-equivariant retraction $r\colon P\to V$ of $i$, we have
		\[
		\mathrm{res}_H^G(\gamma)=\tr_V(r\Phi)\mathrm{res}_H^G(B)\quad \text{in $H^2(H,k)$.}
		\]
	\end{lemma}

\begin{proof}
	Applying $\Hom_k(V,-)$ to \eqref{eq:detector} gives an exact sequence of $G$-modules
	\[
	0\longrightarrow\End_k(V)\xlongrightarrow{i_*}
	\Hom_k(V,P)\xlongrightarrow{\pi_*}\End_k(V)
	\longrightarrow0,
	\]
	and hence an exact sequence
	\[
	H^2(G,\End_k(V))\xrightarrow{i_*}
	H^2(G,\Hom_k(V,P))
	\xrightarrow{\pi_*}H^2(G,\End_k(V)).
	\]
	Define $\widetilde y\coloneqq\Phi_*B\in H^2(G,\Hom_k(V,P))$. Since $\pi\Phi=d$, the class $\widetilde y$ satisfies $\pi_*\widetilde y=d_*B=0$, and so there exists $\Gamma\in H^2(G,\End_k(V))$ such that $i_*\Gamma=\widetilde y$. We set $\gamma\coloneqq\Tr_V(\Gamma)\in H^2(G,k)$.
	
	For every closed subgroup $H\subset G$ over which
	\eqref{eq:detector} splits and every $H$-equivariant
	retraction $r$ of $i$, restricting to $H$ and applying
	$r_*$ gives
	\[
	\mathrm{res}_H^G(\Gamma)
	=(r\Phi)_*\mathrm{res}_H^G(B)\quad\text{in }H^2(H,\End_k(V)),
	\]
	and taking $\mathrm{Tr}_V$ yields the desired equality.
\end{proof}

	\subsection{A representation of $C_p\times C_p$}
	
	Let $C_p=\langle g\rangle$ be a cyclic group of order $p$. We have an exact sequence
	\[0\longrightarrow \F_p\xlongrightarrow{\nu} \F_p[C_p]\xlongrightarrow{g-1} \F_p[C_p]\xlongrightarrow{\epsilon} \F_p\longrightarrow 0,\]
	where $\nu(1)=N=1+g+\dots+g^{p-1}$ and $\epsilon$ is the augmentation map.
	Letting $I\coloneqq (g-1)\subset \F_p[C_p]$ be the augmentation ideal, we obtain two short exact sequences
	\begin{equation}\label{eq:two-extensions}
		\begin{gathered}
			0\longrightarrow\F_p\xlongrightarrow{\nu}\F_p[C_p]\xlongrightarrow{g-1}I\longrightarrow0,\\
			0\longrightarrow I\longrightarrow \F_p[C_p]\xlongrightarrow{\epsilon}\F_p\longrightarrow 0.
		\end{gathered}
	\end{equation}
	Set $W\coloneqq\F_p\oplus I$ and $V\coloneqq \F_p[C_p]\oplus \F_p[C_p]$. Taking the direct sum of the two sequences in \eqref{eq:two-extensions} gives an extension
	\begin{equation}\label{eq:boundary}
		0\longrightarrow W\xlongrightarrow{j}V\xlongrightarrow{q}W\longrightarrow0.
	\end{equation}
	Let $d\coloneqq jq\in \mathrm{End}_{\F_p}(V)$, and define $h\in \End_{\F_p}(V)$ by $h(x,y)=(x,0)$. Then
	\begin{equation}\label{eq:grading}
		d^2=0,\qquad \ker d=\operatorname{im}d=j(W),\qquad dh+hd=d.
	\end{equation}	
	Let $C_p\times C_p=\langle g_1,g_2\rangle$ act on $V$ by $g_1\mapsto g$ and $g_2\mapsto1+d$. The $p+1$ order-$p$ subgroups of $C_p\times C_p$ are parametrized by $\mathbb P^1(\F_p)=\F_p\cup\{\infty\}$, and for every $\lambda\in\F_p$ we set $L_\lambda\coloneqq\langle g_1g_2^\lambda\rangle$, and we set $L_\infty\coloneqq\langle g_2\rangle$.
	
	\begin{lemma}\label{lem:lines}
		For every $\lambda\in\F_p$, the restriction $V|_{L_\lambda}$ is free of rank $2$ over $\F_p[L_\lambda]$. Moreover, $V$ is free of rank $p$ over $\F_p[z]/(z^2)$, where $z$ acts on $V$ by $d$.
	\end{lemma}

	\begin{proof}
		The second assertion follows from $\ker d=\operatorname{im}d$, which holds by \eqref{eq:grading}. Now let $\lambda\in\F_p$, let $\iota_\lambda\colon C\xrightarrow{\sim}L_\lambda$ be the isomorphism sending $g$ to $g_1g_2^\lambda$, and write $V_\lambda$ for the $C_p$-module obtained by pulling back $V$ along $\iota_\lambda$, so that $g$ acts on $V_\lambda$ as $g(1+\lambda d)$. 
		
		We have $d(x,y)=(Ny,(g-1)x)$ for all $x,y\in V$. Since $d^2=0$ and $d$ commutes with $g$, it follows that the action of $g-1$ on $V_\lambda$ is given by \[T_\lambda\coloneqq (g-1)+\lambda gd\in \End_{\F_p}(V),\] and by the binomial theorem, for every $(x,y)\in V$ we have
		\[
		\begin{aligned}
			T_\lambda^{p-1}(x,y)
			&=((g-1)^{p-1}-\lambda(g-1)^{p-2}gd)(x,y)\\
			&=(Nx,Ny)-\lambda(0,Nx)\\
			&=(Nx,N(y-\lambda x)).
		\end{aligned}
		\]
		As $N\cdot \F_p[C_p]=\F_p\cdot N$, we deduce that $T_\lambda^{p-1}$ has dimension two. Moreover, for every $(x,y)\in V$ we have
		\[T_\lambda^p=(g-1)^p+\lambda^p g^p d^p=0.\]
		Thus $T_\lambda$ has exactly two Jordan blocks, both of size $p$, and hence $V_\lambda$ is free of rank two over $\F_p[C_p]$.
	\end{proof}

	\begin{lemma}\label{lem:local-traces}
		There exist an extension of $(C_p\times C_p)$-representations
		\begin{equation}\label{eq:relative-extension}
			0\longrightarrow V\xlongrightarrow{i}P\xlongrightarrow{\pi}V\longrightarrow0
		\end{equation}
		over $\F_p$ and a $(C_p\times C_p)$-equivariant map $\Phi\colon V\to P$ with $\pi\Phi=d$ such that
		\begin{enumerate}
			\item for every $\lambda\in \mathbb P^1(\F_p)$, the restriction of \eqref{eq:relative-extension} to $L_\lambda$ splits, and
			\item for every $L_\lambda$-equivariant retraction $r\colon P\to V$ of $i$, we have
			\[\tr_V(r\Phi)=0\quad(\lambda\in\F_p),\qquad \tr_V(r\Phi)=2\quad (\lambda=\infty).\]
		\end{enumerate}  
	\end{lemma}

	\begin{proof}
		Let $R\coloneqq\F_p[C_p]$, and set $P\coloneqq R[z]/(z^2)\otimes_R V=V\oplus (z\otimes V)$, where $V=R\oplus R$ is regarded as an $R$-module, and let $C_p\times C_p$ act on $P$ by $g_1\mapsto g$ and $g_2\mapsto1+z$. Let $S\colon V\to V$ be given by $S(x,y)=(x,-y)$, and define 
		\[
		\begin{aligned}
			i\colon V&\longrightarrow P,
			&v&\longmapsto z\otimes S(v)-1\otimes dS(v),\\
			\pi\colon P&\longrightarrow V,
			&1\otimes v&\longmapsto v,\quad
			z\otimes v\longmapsto d(v).
		\end{aligned}
		\]
		Since $Sd=-dS$, the maps $i$ and $\pi$ are $(C_p\times C_p)$-equivariant and give an exact sequence \eqref{eq:relative-extension}.
		
		 Define an $R$-linear map
		\[
			\Phi\colon V\longrightarrow P,\qquad v\longmapsto z\otimes h(v)+1\otimes h(dv).
		\]
		Since $d^2=0$, we have $z\Phi(v)=z\otimes h(dv)=\Phi(dv)$, which shows that the map $\Phi$ is $(C_p\times C_p)$-equivariant. By \eqref{eq:grading}, for every $v\in V$ we have
		\[\pi\Phi(v)=\pi(z\otimes h(v)+1\otimes h(dv))=dh(v)+hd(v)=d(v),\]
		that is, $\pi\Phi=d$. Moreover, letting $h_W\in\End_{\F_p}(W)$ be the idempotent given by
		$h_W(a,x)=(a,0)$ for every $(a,x)\in W=\F_p\oplus I$, we have, for every $w\in W$, the identities $dj(w)=0$, $hj(w)=jh_W(w)$, $Sjh_W(w)=jh_W(w)$, and hence
		\[
		\begin{aligned}
			\Phi(j(w))
			&=z\otimes hj(w)
			=z\otimes jh_W(w)\\
			&=z\otimes Sjh_W(w)-1\otimes dSjh_W(w)
			=ijh_W(w),
		\end{aligned}
		\]
		that is, $\Phi j=ijh_W$.
		
		Fix $\lambda\in\mathbb P^1(\F_p)$. By Lemma~\ref{lem:lines}, if $\lambda\neq\infty$ then $V|_{L_\lambda}$ is free over $\F_p[L_\lambda]$, so \eqref{eq:relative-extension} splits over $L_\lambda$, and  if $\lambda=\infty$ then $V$ is free over $\F_p[z]/(z^2)$, so \eqref{eq:relative-extension} splits over $L_\infty$ as well.
		
		Let $r\colon P\to V$ be an $L_\lambda$-equivariant retraction of $i$. If $\lambda\neq\infty$, then $r\Phi$ is an endomorphism of a free $\F_p[L_\lambda]$-module, and hence its trace is zero, since multiplication by any element of $\F_p[L_\lambda]$ has trace zero over $\F_p$. If $\lambda=\infty$, then $r\Phi$ commutes with $d$, and the equality $\Phi j=ijh_W$ gives $r\Phi j=jh_W$. Since $d=jq$, we also have $q(r\Phi)=h_Wq$, and hence we obtain a commutative diagram with exact rows
		\[
		\begin{tikzcd}
			0 \arrow[r] & W \arrow[r, "j"] \arrow[d, "h_W"'] & V \arrow[r, "q"] \arrow[d, "r\Phi"] & W \arrow[r] \arrow[d, "h_W"] & 0\\
			0 \arrow[r] & W \arrow[r, "j"] & V \arrow[r, "q"] & W \arrow[r] & 0.
		\end{tikzcd}
		\]
		By the additivity of trace in short exact sequences, we obtain the equality
		$\tr_V(r\Phi)=2\tr_W(h_W)=2$, as desired.
	\end{proof}

	\subsection{Construction of the characters}
	
	Let $F$ be a field containing a primitive $p$th root of unity $\zeta$, and let $\beta\colon H^1(F,\F_p)\to H^2(F,\F_p)$ be the Bockstein homomorphism associated with the short exact sequence 
	\[0 \longrightarrow \F_p \longrightarrow \Z/p^2\Z \longrightarrow \F_p \longrightarrow 0.\]
	Explicitly, for every $b\in H^1(F,\F_p)$, letting $\widetilde b\colon G_F\to\{0,\ldots,p-1\}$ be the unique set-theoretic lift of $b$, the class $\beta(b)\in H^2(F,\F_p)$ is represented by the $2$-cocycle $c_b\in Z^2(F,\F_p)$ given by
	\begin{equation}\label{eq:bockstein-cocycle}
	c_b(\sigma,\tau)\coloneqq
	\frac{\widetilde b(\sigma)+\widetilde b(\tau)-\widetilde b(\sigma\tau)}{p}\pmod p.
	\end{equation}
	For every $a\in F^\times$, let $\chi_a\in H^1(F,\F_p)$ be the character defined by
	\[\sigma(a^{1/p})=\zeta^{\chi_a(\sigma)}a^{1/p}\qquad\text{for every }\sigma\in G_F,\]
	where $a^{1/p}$ is a $p$th root of $a$. The Kummer map $a\mapsto\chi_a$ induces an isomorphism $F^\times/F^{\times p}\xrightarrow{\sim} H^1(F,\F_p)$. For every $b\in H^1(F,\F_p)$, we have
	\[\beta(b)=b\cup\chi_\zeta\quad\text{in }H^2(F,\F_p);\]
	see \cite[Lemma~4.2]{GM24}.
	
	\begin{lemma}\label{lem:arithmetic}
		Let $F$ be a global field of characteristic different from $p$ such that $\mu_p\subset F$ and $\mu_{p^2}\not\subset F$. There exist $a,b\in H^1(F,\F_p)$ and a finite place $v_0$ away from $p$ such that $a_v,b_v$ are dependent for every place $v$ and
		\begin{equation}\label{eq:arithmetic}
			a_{v_0}=b_{v_0},\qquad \beta(b_{v_0})\neq 0,\qquad
			\beta(b_v)=0\ \text{if }v\neq v_0\text{ and }a_v\neq 0.
		\end{equation}
	\end{lemma}
	\begin{proof}
		Choose, by Chebotarev's theorem, a place $v_0$ away from $p$ where $\kappa(\zeta)_{v_0}$ is nonzero and unramified and, by weak approximation and Kummer theory, a character $b\colon G_F\to \F_p$ ramified at $v_0$, so that $\beta(b_{v_0})=b_{v_0}\cup \chi_{\zeta}\neq 0$ in $H^2(F,\F_p)$; see \cite[Theorem~7.2.15]{NSW08}.
		
		Let $L/F$ be the cyclic extension defined by $b$, and let $\Sigma$ be a finite set of places containing $v_0$, the places above $p$, and the places where $b$ is ramified. Since $b_{v_0}\cup b_{v_0}=0$, the class $b_{v_0}$ is a norm from $L$ locally at $v_0$, and hence, by weak approximation over $L$ and Kummer theory, there exists $a$ in the image of the corestriction map $\operatorname{Cor}_{L/F}\colon H^1(L,\F_p)\to H^1(F,\F_p)$
		such that $a_{v_0}=b_{v_0}$ and $a_v=0$ for every $v\in\Sigma\setminus\{v_0\}$.
		
		Since $a$ is a norm from $L$, we have $a\cup b=0$. If $v$ is a place not in $\Sigma$, then the character $b_v$ is unramified, and hence it spans its own orthogonal for the cup product, so that $a_v$ and $b_v$ are dependent. By construction, they are also dependent at the places of $\Sigma$ and at the archimedean places (where both characters vanish). Finally, an unramified character $\chi\colon G_F\to \F_p$ factors the unramified quotient of $G_F$, which is procyclic, and hence in particular $\beta(\chi)=0$. All in all, this proves \eqref{eq:arithmetic}.
	\end{proof}
	
	\subsection{Proof of Theorem~\ref{thm:main-odd}}

	\begin{proof}[Proof of Theorem~\ref{thm:main-odd}]
	By Lemma~\ref{lem:induction}(2) and (3), we may assume that $F$ contains a primitive $p$th root of unity $\zeta$. We let $V$ be the $(C_p\times C_p)$-representation of \eqref{eq:boundary}, we fix characters $a,b\colon G_F\to C_p$ as in Lemma~\ref{lem:arithmetic}, and we let $G_F$ act on $V$ and \eqref{eq:relative-extension} via  $\phi\coloneqq(a,b)\colon G_F\to C_p\times C_p$.
		
	We first compute the obstruction to lifting $V$ modulo $p^2$. Consider the free $\Z_p[C_p]$-module $\widetilde V\coloneqq\Z_p[C_p]\oplus\Z_p[C_p]$, and define $\widetilde d\in\End_{\Z_p[C_p]}(\widetilde V)$ by $\widetilde d(x,y)=(\widetilde Ny,(g-1)x)$, where $\widetilde N=1+g+\dots+g^{p-1}\in \Z_p[C_p]$.
		Since $(g-1) \widetilde N=0$, we have $\widetilde d^2=0$, and the reductions of $\widetilde V$ and $\widetilde d$ modulo $p$ are $V$ and $d$. Let $\widetilde a,\widetilde b\colon G_F\to\{0,\ldots,p-1\}$ be the unique lifts of $a$ and $b$, respectively, and define
		\[
		\widetilde\rho\colon G_F\longrightarrow
		\GL(\widetilde V/p^2\widetilde V),\qquad
		\widetilde\rho(\sigma)=g^{\widetilde a(\sigma)}(1+\widetilde d)^{\widetilde b(\sigma)}.
		\]
		Since $(1+\widetilde d)^p=1+p\widetilde d$ in $\End_{\Z_p}(\widetilde V)$, we have
		\[
		\widetilde\rho(\sigma)\widetilde\rho(\tau)\widetilde\rho(\sigma\tau)^{-1}
		=1+p\,c_b(\sigma,\tau)\widetilde d,
		\]
		where $c_b$ is the $2$-cocycle representing $\beta(b)$ defined in \eqref{eq:bockstein-cocycle}. Therefore
		\[o(V)=d_*\beta(b)\quad\text{in }H^2(G_F,\End_{\F_p}(V)).\]
		Suppose by contradiction that $V$ lifts modulo $p^2$, so that $d_*\beta(b)=0$. For every place $v$, since $a_v$ and $b_v$ are linearly dependent, $\phi(G_{F_v})$ is contained in $L_\lambda$ for some $\lambda\in\mathbb P^1(\F_p)$, and moreover $\lambda\in \F_p$ if $a_v\neq 0$. By Lemmas~\ref{lem:criterion} and \ref{lem:local-traces}, there exists $\gamma\in H^2(F,\F_p)$ such that
		\[
		\gamma_v=
		\begin{cases}
			0,&a_v\neq 0,\\
			2\beta(b_v),&a_v=0.
		\end{cases}
		\]
		By \eqref{eq:arithmetic}, the class $\gamma-2\beta(b)$ is trivial at every place other than $v_0$, while its localization at $v_0$ equals $-2\beta(b_{v_0})$, which is non-zero because $\beta(b_{v_0})\neq 0$ and $p$ is odd. This contradicts Brauer reciprocity \cite[Theorem~8.1.17]{NSW08}, and hence $V$ does not lift modulo $p^2$, as desired.
		
		It remains to prove that $V$ lifts to $\Z_p$ everywhere locally. For this, let $v$ be a place of $k$. If $a_v\neq 0$, then $\phi(G_{F_v})\subseteq L_\lambda$ for some $\lambda\in\F_p$, and Lemma~\ref{lem:lines} shows that $V|_{L_\lambda}$ lifts to the free $\Z_p[L_\lambda]$-module of rank two. If $a_v=0$, then $G_{F_v}$ acts on $V$ through $b_v$, and hence by Lemma~\ref{lem:lines} the representation $V|_{G_{F_v}}$ is a direct sum of $p$ copies of the extension 
		\[
		0\longrightarrow\F_p\longrightarrow M_v\longrightarrow\F_p\longrightarrow 0
		\]
		classified by $b_v\in H^1(F_v,\F_p)= H^1(F_v,\mu_p)$. By Kummer theory, the class $b_v$ lifts to $H^1(F_v,\Z_p(1))$, that is, the extension lifts to an extension
		\[
		0\longrightarrow\Z_p(1)\longrightarrow \widetilde M_v\longrightarrow\Z_p\longrightarrow 0,
		\]
		so that $\widetilde M_v^{\oplus p}$ is an integral lift of $M_v^{\oplus p}=V|_{G_{F_v}}$.
	\end{proof}
	
	\begin{rmk}\label{rmk:tamely-ramified}
		The representation in Theorem~\ref{thm:main-odd} can be chosen to be tamely ramified at every finite place. Indeed, in Lemma~\ref{lem:arithmetic}, by weak approximation and Kummer theory we may choose $b$ to also be trivial at every place above $p$, and the construction of $a$ already gives $a_v=0$ at these places $v$. Thus, in the case when $\mu_p\subset F$, the resulting representation $V$ is trivial at every place above $p$ and, since its image is a $p$-group, tamely ramified at all other finite places. In the general case, the extension $F(\mu_p)/F$ is tamely ramified, and hence induction of representations preserves tame ramification, so the conclusion follows.
		
		On the other hand, our methods cannot produce an unramified representation which does not lift. Indeed, if $b$ is unramified, then $\beta(b_v)=0$ for every place $v$ of $k$, and hence $\beta(b)=0$ by Brauer reciprocity, so that $o(V)=d_*\beta(b)=0$.
	\end{rmk}
	
	\section*{Acknowledgments}
	
	I thank Olivier Wittenberg for helpful conversations and email exchanges on the Poitou--Tate pairing. I owe Remark~\ref{rmk:tamely-ramified} to a question of Matteo Tamiozzo. I thank Alexander Merkurjev for asking me about the liftability of $Q_8$ and for his interest in the present work.
	
	The proofs of Theorems~\ref{thm:main-even} and \ref{thm:main-odd} have been obtained without AI assistance, and so has the current proof of Theorem~\ref{thm:quaternion}. My original proof of Theorem~\ref{thm:quaternion} was completely different from the one presented here and relied on a quiver-theoretic approach to the representations of generalized quaternion groups. Most of that proof was obtained without AI assistance, but I used AI tools to help carry out a particularly difficult and long computation. I subsequently realized that Theorem~\ref{thm:quaternion} could also be proved using the circle of ideas developed in this paper, and therefore decided to include it here.
	
	This manuscript has been written by myself, but after the writing was complete I have used AI to proofread the manuscript and to improve the exposition.

\end{document}